\UseRawInputEncoding

\newenvironment{dedication}
{
\itshape
\raggedleft
}
{\par
\vspace{3mm}
}

\documentclass{amsart}
\usepackage{amsmath,amssymb}

\usepackage{epsfig}
\usepackage{url}
\usepackage{amssymb}
\usepackage{amsmath}
\usepackage{amsfonts}

\def\Dbar{\leavevmode\lower.6ex\hbox to 0pt{\hskip-.23ex \accent"16\hss}D}

\def\bZ{{\mbox{\bf Z}}}

\begin{document}

\title{Hadamard matrices: skew of orders 276,292 and symmetric 372}
\author {Dragomir {\v{Z}}. {\Dbar}okovi{\'c}}
\address{University of Waterloo, 
Department of Pure Mathematics,
Waterloo, Ontario, N2L 3G1, Canada}
\email{dragomir@rogers.com}
\date{}

\begin{abstract}
The smallest integer $v>0$ for which no skew-Hadamard matrix of order $4v$ is known is $v=69$.
We show how to construct several such matrices. 
We also construct the first examples of a skew Hadamard matrix of order 292 and symmetric Hadamard matrices of order 372.
\end{abstract}

\maketitle

\begin{dedication}
	In memory of my son Dejan Djokovic (1962-2022).
\end{dedication}

\section{Introduction}

According to Table 7.1 of the survey paper of Seberry and Yamada \cite{SY:1992}, published in 1992,
there were only six odd integers $v<100$ for which no skew-Hadamard matrix of order $4v$
was known at that time, namely the integers
$$
47,59,69,81,89,97.
$$
Subsequently, the skew-Hadamard matrices of order $4v$ were constructed in 1994
for $v=81$ \cite{Djokovic:AustralasJC:1994},
in 2004 for $v=59$ \cite{FKS:2004}, and in 2008 for $v=47,97$ \cite{Djokovic:skewHM:188:388:2008}.
In this note we construct 19 skew-Hadamard matrices of order $276$ $(=4\cdot 69)$
and 4 symmetric Hadamard matrices of order $372$ $(=4\cdot 93)$.

Let us make two remarks about the skew-Hadamard matrices. 
First, the case $v=63$ is listed as unknown in the handbook \cite[Table 1.51, p. 277]{CK-Had:2007} published in 2007.
However the existence of a skew-Hadamard matrix of order 
$4\cdot 63$ was known since 1969, as it belongs to
an infinite series of such matrices constructed by Szekeres \cite{Szekeres:1969}.
As this handbook does not list $v=69$ as unknown, it is probable that this was just a misprint: 63 should be replaced by 69?
Second, in the more recent book \cite[Table 9.2, pp. 198-200]{SY:2020}, the cases $v=39,49,65$ are listed as unknown.
However the corresponding skew-Hadamard matrices have been constructed long ago in \cite{Djokovic:PEF:1992}. 

\section{The first skew-Hadamard matrices of order 276 and 292}

As far as we know, the smallest integer $v>0$ for which no skew-Hadamard matrix of order $4v$ is known is
$v=69$ \cite[p. 1436)]{BDKG:2019}.
In this section we construct several such matrices and one of
order 292.

Our construction uses the Goethals-Seidel array (GS-array) shown below

\begin{equation*} 
\left[ \begin{array}{cccc}
A_0    & A_1 R & A_2 R & A_3 R \\
-A_1 R & A_0   & -RA_3 & RA_2  \\
-A_2 R & RA_3  & A_0   & -RA_1 \\
-A_3 R & -RA_2 & RA_1  & A_0  \\
\end{array} \right].
\end{equation*} 

We shall assume that the $A_i$ are circulants and $R$ is the 
back-diagonal identity natrix (i.e. the matrix obtained from the
identity matrix by reversing the order of rows). The circulants 
are obtained from a cyclic difference family $\{X_0,X_1,X_2,X_3\}$ with parameters
$$
(v=69; k_0=34, k_1=34, k_2=31, k_3=27; \lambda=57).
$$
For instance, for the first row $(a_0,a_1,\ldots,a_{v-1})$ of $A_0$ we have $a_i=-1$ if $i\in X_0$ and $a_i=1$ otherwise.
Moreover it is required that the block $X_0$ is skew, i.e. $a_0=1$ and $a_i + a_{v-i}=0$ for $i=1,2,\ldots,34$.

A special feature of our difference families is that they  
break up into two pieces $\{X_0,X_1\}$ and $\{X_2,X_3\}$ which are also difference families.

First, we need a difference family $\{X_0,X_1\}$ with parameters $(69; 34,34; 33)$ with $X_0$ skew.
This is provided by the well known family of Szekeres difference sets \cite{Szekeres:1969,SY:2020}:
\begin{eqnarray*}
	X_0 &=& \{1,2,6,7,9,13,14.16,17,18,21,27,31,34,36,37,39,40,    \\
          & &     41,43,44,45,46,47,49,50,54,57,58,59,61,64,65,66\};   \\
        X_1 &=& \{1,4,5,7,9,10,11,12,15,17,18,19,24,26,27,28,30,39,    \\
            & &   41,42,43,45,50,51,52,54,57,58,59,60,62,64,65,68\}.  
\end{eqnarray*}
Note that $X_0$ is skew and $X_1$ is symmetric.
Further, all 19 difference families share the same first two blocks, $X_0$ and $X_1$.

Second, we need a difference family $\{X_2,X_3\}$ with parameters $(69; 31,27; 24)$. 
This is exactly the parameter set for a D-optimal design of order $2\cdot 69=138$. 
In a joint paper with I. Kotsireas \cite[Section 4.2]{DK:D-optimal:2015}, 
we have constructed 19 nonequivalent such difference families.
Anyone of them can be used in our construction. As an example, let us choose the first one:
\begin{eqnarray*}
	X_2 &=& \{0,1,3,4,6,9,10,11,13,14,17,18,20,22,26,28,29, \\
	     &&   32,33,34,39,41,43,45,46,48,51,59,60,62,63\},  \\
	X_3 &=& \{0,2,3,4,8,9,10,11,12,15,16,17,21,25,26, \\
	     &&   32,33,35,36,37,39,41,46,51,54,57,59\}. 
\end{eqnarray*}
By constructing the circulants $A_i$ from the blocks $X_i$ and
by plugging the $A_i$ into the GS-array we obtain a skew-Hadamard
matrix of order 276.

Consequently, the smallest positive integer $v$ for which the existence of a skew-Hadamard 
matrix of order $4v$ is still undecided is now $89$.

For the readers convenience, we provide (for the difference family chosen above)
the first rows of the blocks $A_i$:
\begin{eqnarray*}
& +  - - + + + - - + - +  + + - - + - - - + +  - + + + + + - + + +  - + + -  \\
&    + - - + - - - + - -  - - - + - - + + + -  + + - - - + - + + -  - - + +; \\
& +  - + + - - + - + - -  - - + + - + - - - +  + + + - + - - - + -  + + + +  \\
&    + + + + - + - - - +  - + + + + - - - + -  + + - - - - + - + -  - + + -; \\
& -  - + - - + - + + - -  - + - - + + - - + -  + - + + + - + - - +  + - - -  \\
&    + + + + - + - + - +  - - + - + + - + + +  + + + + - - + - - +  + + + +; \\
& -  + - - - + + + - - -  - - + + - - - + + +  - + + + - - + + + +  + - - +  \\
&    - - - + - + - + + +  + - + + + + - + + -  + + - + - + + + + +  + + + +. \\
\end{eqnarray*}
(The $+$ and $-$ signs stand for $+1$ and $-1$, respectively.)

As far as we know, the odd integers $v>0$ less than 200 for which the existence
of skew-Hadamard matrices of order $4v$ is still undecided are the following:
$$
89,  101, 107, 119, 149, 153, 167, 177, 179, 191, 193.
$$
After taking into account the papers \cite{Djokovic:skewHM:188:388:2008,Djokovic:JCD:2008}
(and correcting the hypothetical misprint mentioned earlier),
this list agrees with \cite[Table 1.51, p. 277]{CK-Had:2007}.

However the above list is still questionable. Indeed
there are many places in mathematical literature where it is
asserted that skew Hadamard matrices of order 292 are known.
Dmitrii Pasechnik pointed out that none of the places that he consulted explains how
to construct such a matrix. As my search also failed, it seems likely that no such matrix
has been known prior to this year (2024).
See also the preprint \cite{Mat-Dima} of Matteo Cati and D. Pasechnik.

We now give the first example of such a matrix (constructed in january of this year).
Set $v=73$ and let $H=\{ 1,2,4,8,16,32,37,55,64 \}$, the subgroup of order 9 of $\bZ_v^\ast$.
The four blocks
\begin{eqnarray*}
	& X_0=H\cdot\{5,9,11,25\}, & X_1=H\cdot\{11,13,17,25\},\\
	& X_2=H\cdot\{5,9,13,17\}, & X_3=H\cdot\{0,1,3,13\}
\end{eqnarray*}
form a difference family in $\bZ_v$ having the parameters
$(73;36,36,36,28;63)$. Since $X_0$ is of skew type, the
corresponding cyclic matrices $A_i$ can be plugged into
the GS-array to obtain the desired skew Hadamard matrix.

\section{The first examples of symmetric Hadamard matrices of order 372}

There are only three odd integers $v<100$ for which no 
symmetric Hadamard matrix of order $4v$ is known \cite{Djokovic:ICD:2022}.
These integers are $65,89$ and $93$. We have constructed four 
non-equivalent such matrices for $v=93$.
Our construction uses the so called {\em propus array} shown below

\begin{equation*} 
\left[ \begin{array}{cccc}
-A_0    & A_1 R & A_2 R & A_3 R \\
 A_2 R  & RA_3  & A_0   & -RA_1 \\
 A_1 R  & A_0   & -RA_3 & RA_2  \\
 A_3 R  & -RA_2 & RA_1  & A_0   
\end{array} \right].
\end{equation*} 

Note that this array is obtained from the GS-array by multiplying the first column
by $-1$ and switching the second and third rows.
If the blocks $A_i$ are circulant matrices and 
$A_1=A_2$ then the propus array gives a symmetric matrix.

In our four examples the blocks $A_i$ are circulants of size 
$v=93$. Since $A_1=A_2$ we need only to provide the first 
rows of the blocks $A_0,A_1$ and $A_3$. We label the $93$ positions of the first
row by integers $0,1,2,\ldots,92$ in that order.
It suffices to list the positions of the entries $-1$ in 
the first row. These positions are provided by a 
difference family $(X_0,X_1,X_2,X_3)$ in $\bZ_{93}$ with parameters $(93; 45,41,41,41; 75)$ and with $X_1=X_2$.
In the families that we use each block $X_i$ is a union of orbits of the subgroup $H=\{1,25,67\}$ of $\bZ_{93}^*$.
Instead of listing the elements of $X_i$ we list below just 
a set of representatives of the $H$-orbits contained in $X_i$ for $i=0,1,3$.
Our four non-equivalent difference families are given by the following four listings:

\begin{eqnarray*}
&[[0,1,3,4,9,11,16,17,18,20,26,31,40,44,48,55,62], \\
&~[2,3,9,13,16,20,22,29,31,36,43,44,48,51,62], \\
&~[1,2,3,8,10,13,17,26,31,36,37,40,43,51,62]]; \\
& \\
&[[0,6,8,12,13,16,17,22,24,31,33,36,40,44,47,51,62], \\
&~ [1,4,5,6,9,11,13,20,31,33,44,48,51,55,62], \\
&~ [1,3,5,6,8,11,17,18,22,26,31,37,40,55,62]]; \\
& \\
&[[0,2,6,8,12,17,20,22,24,31,33,36,40,43,51,55,62], \\
&~ [1,2,3,5,8,9,10,13,16,24,29,31,51,55,62], \\
&~ [1,4,5,6,9,12,13,20,22,31,36,43,47,51,62]]; \\
& \\
&~[[0,1,2,6,9,10,12,16,26,29,31,33,36,43,44,48,62], \\
&~ [2,4,5,6,13,16,29,31,33,36,40,43,44,51,62], \\
&[1,3,6,8,9,12,13,17,20,22,31,43,51,55,62]].
\end{eqnarray*}

For instance the difference family of the first example is explicitly:

\begin{eqnarray*}
&X_0=&\{0,1,3,4,7,9,11,15,16,17,18,20,23,25,26,28,31,35,38,39,\\
&  &40,44,45,48,49,53,54,55,58,62,65,67,68,70,73,75,76,77,78, \\
&  &82,84,86,89,90,92\}, \\
&X_1=&\{2,3,9,13,15,16,20,22,28,29,31,34,35,36,38,39,41,43,44,\\
&  &45,46,48,49,50,51,52,54,62,63,65,66,69,74,75,77,79,83,84, \\
&  &85,87,91\}, \\
&X_2=&X_1, \\
&X_3=&\{1,2,3,8,10,13,14,15,17,19,23,25,26,31,34,36,37,40,41,\\
&  &43,46,50,51,52,53,61,62,63,64,66,67,68,69,70,71,75,76,87,\\
&  &88,91,92\}.
\end{eqnarray*}

In all four difference families the block $X_0$ is symmetric i.e. $-X_0=X_0$ in $\bZ_{93}$
and the representative $0$ represents the trivial $H$-orbit, namely  $\{0\}$.

\section{Acknowledgements}
This research was enabled in part by support provided by SHARCNET (http:// \\
www.sharcnet.ca) and the Digital Research Alliance of Canada (alliancecan.ca).

\end{document}